\documentclass[3p]{elsarticle}
\usepackage{amssymb}
\usepackage{amsmath}
\usepackage{theorem}

{\theorembodyfont{\rmfamily}

}

{\theorembodyfont{\rmfamily}

}
{\theorembodyfont{\rmfamily}

}
{\theorembodyfont{\rmfamily}

}

\begin{document}

\title{A hybrid heuristic algorithm with application to a graphical interface for vehicle routing optimization in an agricultural cooperative}

\author[rg]{Roque M. Guiti\'an de Frutos}
\ead{roqueguitian@hotmail.com}

\author[cm]{Balbina V. Casas-M\'{e}ndez}
\ead{balbina.casas.mendez@usc.es}

\cortext[cor1]{Corresponding author: Balbina Virginia Casas M\'endez. Facultad de Matem\'aticas, R\'ua Lope G\'omez de Marzoa s/n. Campus Vida. Santiago de Compostela (A Coru\~{n}a). Spain. Postal code: 15782. Phone: +34 881813219.}

\address[rg]{Department of Mathematical Analysis, Statistics, and Optimization, University of Santiago de Compostela, Spain.}

\address[cm]{MODESTYA Research Group, Department of Mathematical Analysis, Statistics, and Optimization, University of Santiago de Compostela, Spain.}

\begin{abstract}

This work considers the problem of optimize the routes of the vehicles used by a real agricultural cooperative that distributes animal feed among the partners.
Because solving the exact model is computationally burdensome, we propose to implement heuristic algorithms.
Thus, firstly an initial solution is obtained through an insertion heuristic algorithm. Secondly, we design a simulated annealing algorithm for possible improvements on the initial solution.  Finally, we built a graphical interface to quickly and easily interact with the system. 

\end{abstract}


\maketitle

\bigskip

\section{Introduction}

The study of decision support systems for logistics is an extremely active area of research and applications in modern Operations Research. In this framework, the planning and scheduling of routes for different purposes is a main topic. From a mathematical point of view, such problems are complex and very difficult to solve exactly. Therefore, many heuristic algorithms have been proposed to overcome these difficulties taking advantage of the underlying combinatorial structure of each particular problem. 

This work is motivated by a real problem faced by a livestock feed company. AIRA is an agricultural cooperative located in Taboada (cf. http://www.aira.es/), a town in the Spanish north western region of Galicia. They produce four different kinds of livestock feeds and serve them to their customers. There are 1500 stockbreeders (customers) spread in about 60 municipalities of the rather vast surrounding area as we can see in Figure 1.

\begin{figure}
\begin{center}
	\includegraphics[height=4.5cm]{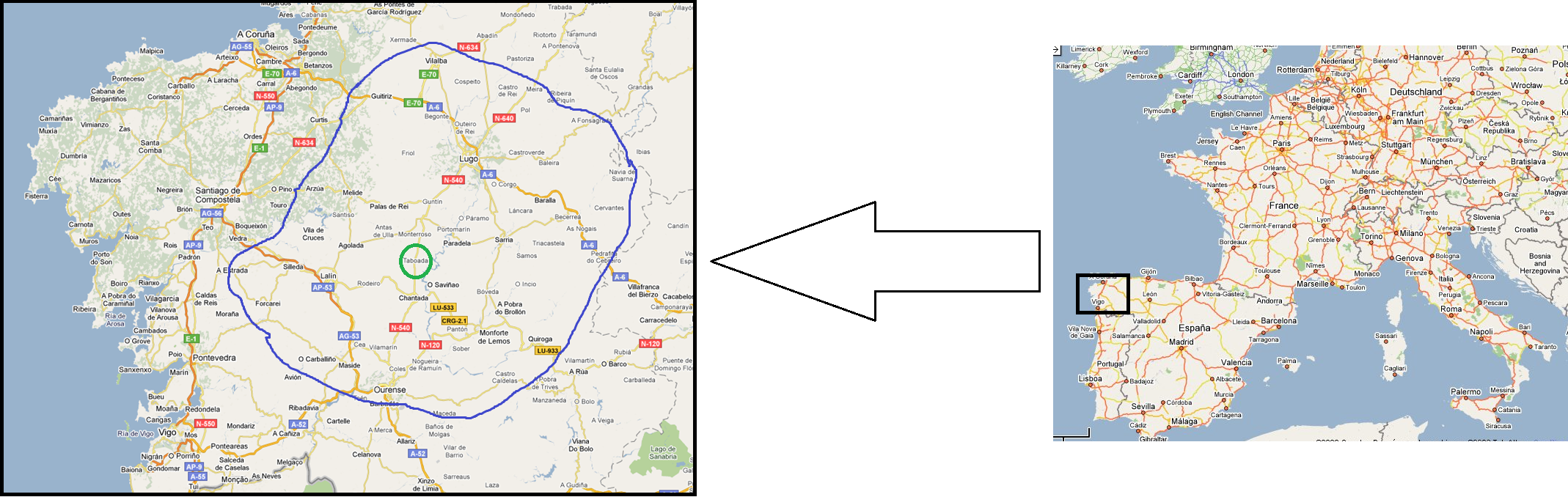}
\end{center}
\caption{Location of agricultural cooperative AIRA and area of influence.}
\end{figure}

The stockbreeders place their orders for the different feeds with the cooperative and each order is either urgent, which means that they need to be served the following day, or has a three or four days window to deliver it. Urgent orders are more expensive for the customers. Usually, stockbreeders order once or twice a month and most of the times, they are only in need of one type of feed. The distribution of feeds is carried out by the fleet of trucks owned by the cooperative and each truck driver is paid depending on the distance covered and the load carried. Additionally, the trucks have different size and are compartmentalized in several hoppers (three to five) which have distinct capacities. Furthermore, given their size, some trucks cannot access some livestock farms. Moreover, for technical reasons, each hopper can only contain one type of feed and cannot be shared among stockbreeders. 

The goal of this study is to provide the cooperative a tool that automatically proposes routes for each truck that satisfying the restrictions, first maximizes the total food delivered per working day and second minimizes the transportation cost of that amount. Besides, the transportation cost paid to each truck driver is a fixed amount for each unloading plus another amount which is proportional to the distance covered and the amount of feed carried on each route. Furthermore, the ratio of this proportion depends on the distance covered.

A research group of Agroforestry Engineering of the University of Santiago de Compostela designed a Global Positioning System (GPS) to monitor the routes of the trucks and the cooperative managers are hoping to improve the whole system in the years to come (cf. Amiama-Ares et al., 2014 \cite{ami}). Indeed, this routing work is a contribution to a bigger project which aims at automating the whole distribution process of the cooperative. 
The GPS that are currently being used provide us all the geographical data we are using to solve this problem, that is to say the distance between all the partners and the cooperative but also the corresponding time necessary to make this route. The other data used here are the quantities of feeds ordered, the capacities of the hoppers, and whether or not a truck can reach a barn. This problem can be formulated using a mathematical programming model\footnote{This note does not include the formulation of the mathematical model of the problem, but the authors can provide it if requested.} that can be written in any language like AMPL (cf. http://ampl.com/). Models prepared in AMPL can be run on the NEOS server (cf. https://neos-server.org/neos/) or using a free or commercial license. 

Ruiz et al. (2004) \cite{rui} is also inspired by the real problem of another Spanish cooperative that manufactures animal feed. Their approach generates all feasible routes by means of an implicit enumeration algorithm, in a first step, and secondly they use an integer programming model that selects the optimal routes. They reduce the initial complexity of the problem by
clustering of customers who have to serve. A recent example that considers the problem of vehicle capacity and different accessibility of different trucks to different customers is found in the work of Derigs et al. (2013) \cite{der}. 

In this scope, we first treat the problem through a nonlinear multiobjective model that for small instances could be solved by a lexicographic procedure. In this way, we confront the task of adapting more classical vehicle routing models to the specific one of AIRA. After that, we can use some algorithm that provides an ``exact solution''. However, we reject in practice this approach since we consider examples with only two trucks and six partners that require more time than three hours when using the NEOS optimization server. In this way, we opted for a hybrid heuristic approach, which first uses an insertion heuristic algorithm to obtain an initial solution (cf. Clarke and Wright, 1964 \cite{cla}, and Mole and Jameson, 1976 \cite{mol}). Then, improvement processes are performed in the initial solution by using a simulated annealing metaheuristic (cf. Kirkpatrick et al., 1983 \cite{kir}, and Nunes de Castro, 2006 \cite{nun}).

\section{Programming model of the problem}

Although we will focus on a methodology based on heuristic algorithms, we briefly explain the basic ideas of the underlying mathematical programming model. This model was presented in \'Alvarez et al. (2010) \cite{alv}. The problem described can be modeled as a modified vehicle routing problem that needs continuous variables representing the load carried by the trucks, it includes time windows associated with agent urgencies, and it incorporates several objectives.

\subsection{Problem data}

In order to model our problem, we use a set of stockbreeders, a set of trucks, and for each truck, its corresponding sets of possible routes and hoppers.
As for the data, we know the capacities of each hopper, the demand of the different customers, and the urgency of their order as well as the distances between stockbreeders and the time it takes to go from one to another. We also know whether a truck is able to reach a specific barn and its maximum authorized charge. Moreover, the cooperative provided us with all the related cost data such as the price of each unloading and the variable cost of the transportation of a certain quantity of feed on a given distance.

\subsection{Decision variables}

The problem as we modeled it involves two different kinds of variables. The first one is
a set of binary variables which equals 1 if there is a route that a truck follows to go from a stockbreeder to another, and 0 otherwise. The second one is a set of continuous variables, in [0,1], which
represents the percentage of a certain hopper that is filled by feed demanded for a certain stockbreeder.

\subsection{Objective functions}

The general objective of this problem is to minimise the cost of food transportation when we try to meet each stockbreeder's order as well as possible.
For doing so, following requirements of the cooperative, we need to divide the problem in two steps: first we maximize the quantity of food carried by the different trucks and then we minimize the cost of the deliveries. In this second step, the term that we minimise is the sum of:
the cost generated by each unloading, which is proportional to the number of stockbreeders visited, and 
the cost generated by the transportation of the feeds which is the sum of a fixed cost proportional to the amount of food carried and a variable cost that depends on the distance covered on the route and again the quantity of food carried on this route. We should note that, actually, the fixed cost of the food transportation is not taken into account in the minimisation as it is considered constant once it has been determined by the previous maximisation.

\subsection{Constraints}

The most common constraint of a vehicle routing problem is the flow condition. Here, every truck that visits a stockbreeder leaves afterwards and all the routes start from the cooperative headquaters. This way, we are sure that if a route really exists, it will end back at the cooperative.
In order to avoid pointless trips, we ask that if a truck goes to a stockbreeder, it carries food for him.
If a truck cannot reach a barn, we make sure that it does not visit the corresponding stockbreeder.
A truck driver cannot work for more than 9 hours a day, so each truck cannot travel more than 9 hours.
A truck cannot cover more than a certain amount of kilometres a day.
The load of each truck on each route is bounded.
To the stockbreeders whose order is urgent (no days left to deliver), we give the quantity they asked for.
To the stockbreeders whose order is less urgent and could still be delivered another day, we give as much as possible.
A hopper cannot be shared among stockbreeders and cannot contain different type of feeds. Thus, each hopper contains only one feed for one stockbreeder.

\subsection{Example}

In order to illustrate the model we elaborated, we consider one small real instance of data. 
In this example, we consider real data provided by the cooperative AIRA. We take 2 similar trucks, with 5 hoppers who can respectively contain up to 3, 3.7, 3.8, 3.7 and 3 tons of feeds. The trucks cannot be carrying more than 11.6 tons. Let's suppose that the trucks can follow one route a day.
Moreover, there are 5 stockbreeders all of them with urgent orders on only one type of feed. The truck is able to deliver food to all the stockbreeders. Besides, to simplify the problem and to speed up its resolution, we suppose that the objective is just minimize the total distance of travelling (also because all orders are urgent). Note that the final number of used hoppers is bigger than the number of served stockbreeders and the solution shows how we have to share each order among different hoppers.
The orders are as shown in Table 1 and the distances between the different barns are as depicted in Table 2.

\begin{table}[h]
\begin{center}
\caption{Orders of the different stockbreeders (in tons of feed).}
{\small
\begin{tabular}{lcc}
\hline
Orders&\hspace*{0.2cm}Feed 1 &\hspace*{0.2cm}Days left to deliver\hspace*{0.2cm}\\
\hline
Stockbreeder 1&\hspace*{0.2cm}3.300&\hspace*{0.2cm}0\\
Stockbreeder 2&\hspace*{0.2cm}2.951&\hspace*{0.2cm}0\\
Stockbreeder 3&\hspace*{0.2cm}3.003&\hspace*{0.2cm}0\\
Stockbreeder 4&\hspace*{0.2cm}3.016&\hspace*{0.2cm}0\\
Stockbreeder 5&\hspace*{0.2cm}2.496&\hspace*{0.2cm}0\\
\hline
\end{tabular}
} 
\label{Orders}
\end{center}
\end{table}

\begin{table}[h]
\begin{center}
\caption{Traveling distances between pairs of stockbreeders (in kilometers).}
{\small
\begin{tabular}{lccccccccc}
\hline
Min&Coop \hspace*{0.00cm}&Stockb. 1\hspace*{0.00cm}&Stockb. 2\hspace*{0.00cm}&Stockb. 3\hspace*{0.00cm}& Stockb. 4\hspace*{0.00cm} & Stockb. 5\hspace*{0.00cm}\\
\hline
Coop.    &0\hspace*{0.00cm}&28\hspace*{0.00cm}&69\hspace*{0.00cm} &64\hspace*{0.00cm} &27                   &17\\
Stockb. 1&\hspace*{0.00cm} &0 \hspace*{0.00cm}&67 \hspace*{0.00cm}&62 \hspace*{0.00cm}&20 \hspace*{0.00cm}  &20\\
Stockb. 2&\hspace*{0.00cm} &\hspace*{0.00cm}  &0\hspace*{0.00cm}  &7  \hspace*{0.00cm}&74 \hspace*{0.00cm}  &58\\
Stockb. 3&\hspace*{0.00cm} &\hspace*{0.00cm}  &\hspace*{0.00cm}   &0\hspace*{0.00cm}  &69 \hspace*{0.00cm}  & 53\\
Stockb. 4&\hspace*{0.00cm} & \hspace*{0.00cm} &\hspace*{0.00cm}   &\hspace*{0.00cm}   &0\hspace*{0.00cm}    & 25\\
Stockb. 5&\hspace*{0.00cm} & \hspace*{0.00cm} &\hspace*{0.00cm}   &\hspace*{0.00cm}   &\hspace*{0.00cm}     &0\\
\hline
\end{tabular}
} 
\label{Traveling times}
\end{center}
\end{table}

The implementation of the model was done in the language AMPL\footnote{Detailed implementation with AMPL can be obtained from authors upon request.} (cf. Fourer et al., 2003 \cite{fou}). We solved the problem with a free version of the solver KNITRO 6.0.0 (cf. Byrd et al., 2006 \cite{byr}), which uses a version of the interior point algorithm for non linear problems, obtained in the NEOS optimization server. When dealing with integer variables, the solver uses a branch and bound algorithm which can only give a locally optimal solution.

The locally optimal solution the solvers gives us after 109.91818 secs is:

\begin{itemize}
	\item The truck 1 goes from Taboada to stockbreeder 5, then to stockbreeder 3, and then to stockbreeder 2 and carries:
\begin{itemize}
	\item 1.475499 tons of feed to the stockbreeder 2 in hooper 1,
	\item 2.496 tons of feed to the stockbreeder 5 in hooper 2,
	\item 1.5505862 tons of feed to the stockbreeder 3 in hooper 3,
	\item 1.4524165 tons of feed to the stockbreeder 3 in hooper 4, and
		\item 1.475499 tons of feed to the stockbreeder 2 in hooper 5.
\end{itemize}
\item The truck 2 goes from Taboada to stockbreeder 4 and then to stockbreeder 1 and carries:
\begin{itemize}
	\item 1.508001 tons of feed to the stockbreeder 4 in hooper 1,
	\item 3.2999998 tons of feed to the stockbreeder 1 in hooper 3, and
	\item 1.508001 tons of feed to the stockbreeder 4 in hooper 5.
\end{itemize}\end{itemize}

In this configuration the total distance of the deliveries is 221 kilometers.

\section{Routing implementation through heuristic algorithms}

Because our problem is hard to solve exactly, we design and apply a hybrid heuristic approach that combines an insertion algorithm to obtain an initial solution and after a simulated annealing metaheuristic to improve it because these kind of strategies have shown to be effective tools for solving hard combinatorial optimization problems.

\subsection{Heuristic algorithm to obtain an initial solution}

It aims to generate an initial solution to the problem, although not necessarily optimal. It itself is a quality solution with which to work in future optimizations. The package would have attributes that implement the data of the problem. The algorithm works on each journey building until all orders are planned. This is the pseudocode defined below. 

1. The programming of each journey begins with the choice of a customer without service (seed client) and an initial truck (seed truck). If there are no customers to serve, the algorithm terminates. If no truck can serve any customer we assume that all trucks are making their journey completely assembled and then we starts the new schedule (with empty vehicles) the next day to serve customers outstanding. If 365 days of programming are completed and there are still customers without serving, the algorithm ends. 2. It is constructed the way that parts of the central repository, visits the seed customer and returns to the company. For the outward journey, the hoppers,
of the truck chosen in the previous step, are loaded with the order data in the manner
more efficient. 3. It analyzes all customers who have not been served and insertion possibilities they have in each subpath of the path. The client and position are obtained with less cost of insertion. 4. It is inserted into the journey the customer with the  obtained position. Hoppers are loaded in all subpaths that precedes this position. Let's suppose it is found that any restriction of the problem is breached, if so return to step 3 discarding this client, otherwise return to step 2 to continue inserting clients. 5. If in step 2 any customer can be inserted, consider the path defined. We return to step 1 to define more routes for more trucks. In the case of loading a truck hopper, if an order can not fit in the vehicle, it is divided into two, one is loaded into the free space of the truck and is considered served and the other not. 

The user can select the seed customer choice strategy among those shown below.

1. Customer farthest. 2. Customer who has more orders that have not yet been addressed.
3. Any client. Sometimes, random choice produces solutions of surprising quality.

The user could also select the strategy of truck choice used in each case, among which are shown below.

1. The truck with lower mileage. 2. The higher capacity truck. 3. Any truck. Sometimes, random choice produces amazing quality solutions again.

Figure 2 resumes the algorithm.

\begin{figure}
\begin{center}
	\includegraphics[height=8.5cm]{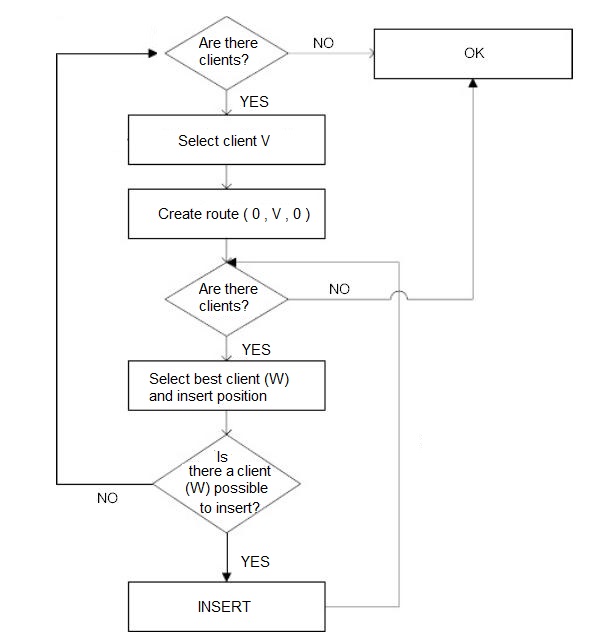}
\end{center}
\caption{The insertion algorithm.}
\end{figure}

\subsection{The simulated annealing algorithm}

At each iteration, the search process moves
of the current solution to a neighboring solution. Let $T$ a parameter that measures the tendency to accept the current candidate. It is called temperature. Among the immediate neighbors of the current solution, select one of
randomly to become future proof solution. If the candidate under consideration is better than the current solution
is accepted. If the candidate under consideration is worse, the chance of acceptance will depend on how much worse it is and on the temperature value, which will evolve over the execution of the algorithm (it will decrease and thus the associated probability). For more details about simulated annealing algorithms see for example Nunes de Castro (2006) \cite{nun}. 

In order to implement an efficient software, seven ways to obtain neighbors solutions are designed and programmed that are  randomly interspersed. Only in this way, we will be able exploring a sufficiently large part of its region of potential solutions
where the global optimum can be located if the execution time permits. One procedure is as follows. 1. Getting a new neighbor moving the orders of a day, truck, customer, and journey randomly selected to another different truck, journey and day
randomly selected to. 

The remaining six procedures have a similar flavour.
2. Getting a new neighbor moving orders of one day, truck, journey, and client randomly selected to a different path of the same truck and day. 3. Getting a new neighbor moving the orders of a day, truck, customer, and journey to a truck and random journey of the same day. 4. Getting a new neighbor exchanging the orders of a day, truck, journey, and random customer by another random customer and different path
of the same truck and day. 5. Getting a new neighbor exchanging orders of a day, truck, customer, and journey randomly selected by other different day, truck, customer, and journey randomly selected. 6. Getting a new neighbor choosing randomly one day, a truck
scheduled for that day, a path of that truck and two clients served
in that way. We exchanged the order of those customers for
that truck, this way and that day. 7. Getting a new neighbor choosing randomly one day, two trucks
scheduled for that day, a path for each truck and a client
served in each way. We exchange the served farmers of a trip
and truck to other.

Figure 3 resumes the algorithm.

\begin{figure}
\begin{center}
	\includegraphics[height=8.5cm]{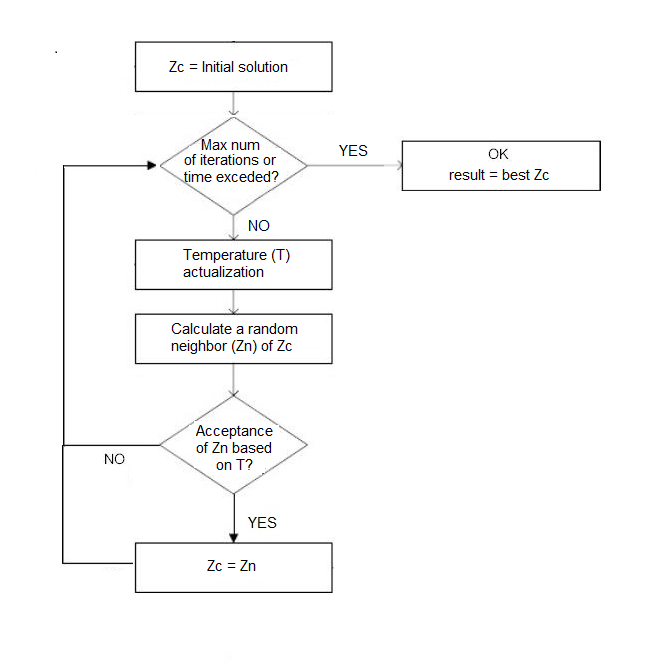}
\end{center}
\caption{The simulated annealing algorithm.}
\end{figure}

\subsection{Graphical interface development}

An application is built making use of the programming language JAVA that includes the heuristic algorithms presented above. With this graphical interface, the user
can introduce data, see results, and in general interact with the system quickly and easily. For this, it has been used
the graphic developer included in the NetBeans development environment (cf. https://netbeans.org/).
The volume of code programmed
for this is about 5818 instructions.

Within the usability criteria for the system implementation will be sought to encourage the expeditious and effective interaction with the user, minimizing
delays and errors arising from the interface.

\begin{figure}
\begin{center}
	\includegraphics[height=8.5cm]{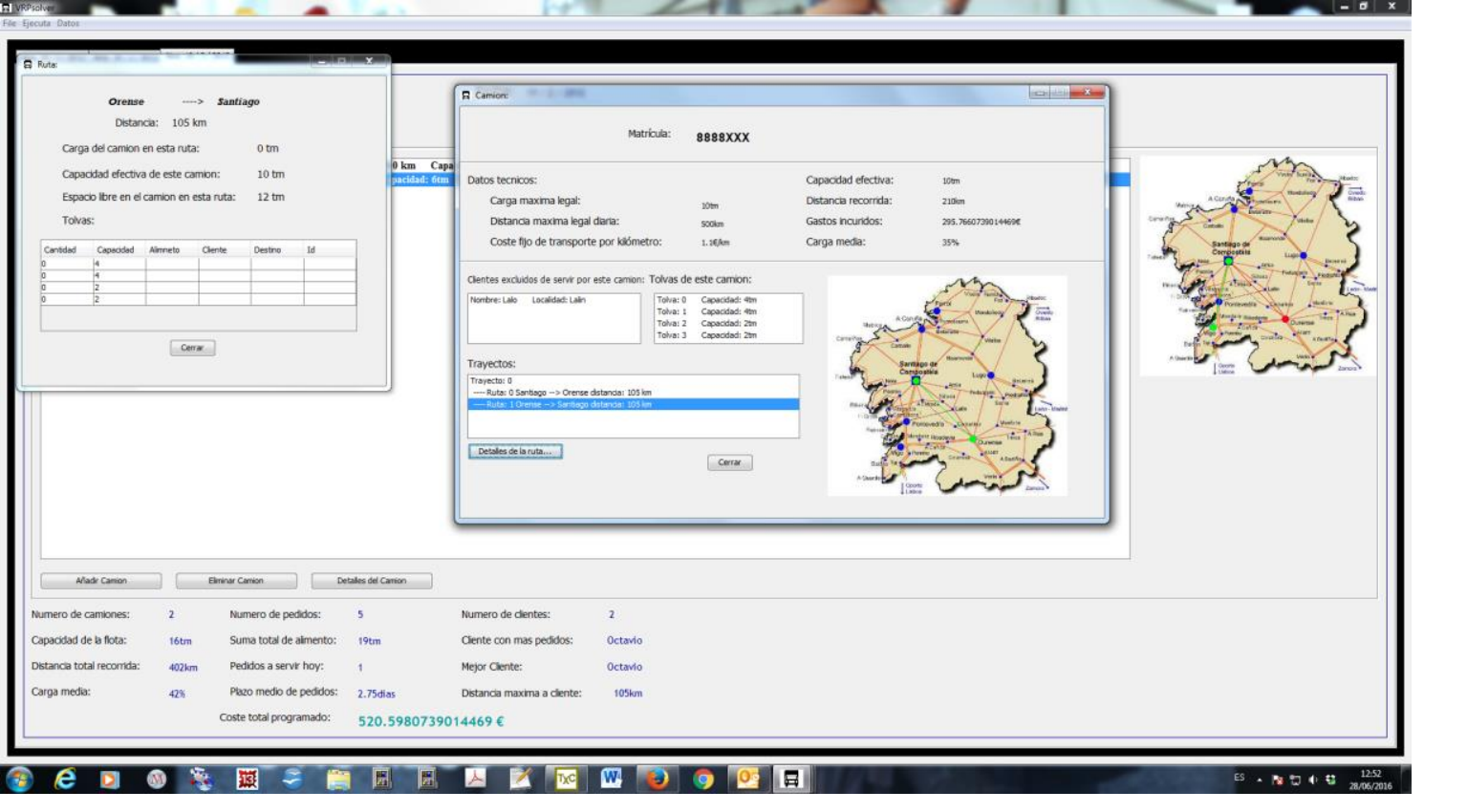}
\end{center}
\caption{Appearance of the graphical interface.}
\end{figure}

Figure 4 below (in Spanish) shows the application after running the algorithms. Selected a truck, we can see on the map the route it must follow. We can see in detail the data from one of the trucks, the drawing of the path, and the distribution of food in the different hoppers.

\section{Some experimental results}

This section describes results obtained from real examples and it presents evidence made to verify the requirements of our application.

\subsection{Former methodology vs. here proposed methodology}

In order to illustrate the graphical interface we elaborated and the algorithms we implemented, we consider an  instance of real data. We plan the work of the cooperative trucks 
for a day of work.

We chose to work with a fleet of 2 trucks, a little one with 5 hoppers who can respectively contain up to 3000, 3700, 3800, 3700, and 3000 kilograms of feed and a bigger one with also 5 hoppers whose capacities are 4000, 3000, 1700, 4500, and 3000 kilograms of feed. Both trucks cannot be carrying more than 15300 kilograms at the same time. Moreover, there are  17 stockbreeders located in different towns. We assume that both trucks are able to deliver food to all partners. The orders are known. We also consider the distance among the different barns. We fix the maximum distance covered by a truck daily, the cost of an unloading, the fixed cost of food transportation, and the variable cost of food transportation.

We introduce the data in the application and with all the information we create a new case. We carry out the execution of both algorithms. With respect to performance parameters of the insertion algorithm for calculating the initial solution, seed customer is taken as one having more orders not yet addressed and the truck to be used is chosen randomly. As for the parameters of simulated annealing algorithm execution, the maximum number of iterations is 750000 and the maximum waiting time is 5 minutes. 

The provided result indicates that the total distance to be traveled by the two trucks is 339 kilometers with an associated cost of 499.21 euros. The distance actually traveled by the trucks was 378.61 kilometers representing an increase of 11.68 percent.

\subsection{Software quality and efficiency of algorithms}

With the numerical results of different tests, they have been constructed several graphs that relate the variables and parameters of the executed algorithms. From
their analysis, we can obtain information about the efficiency of the application.
Figure 5 shows how the system behaves as it faces problems of different sizes. Thus, the vertical axis indicates the percentage of
improvement between the initial solution and the optimized solution, while the 
horizontal axis shows the number of problem data. The optimization runs for
a time constant of 30 seconds.

\begin{figure}
\begin{center}
	\includegraphics[height=8cm]{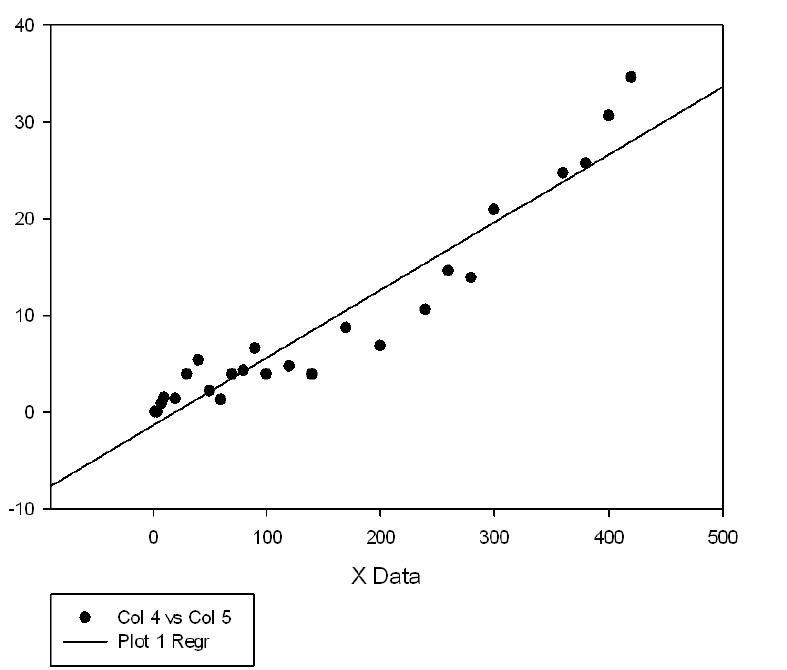}
\end{center}
\caption{The system behavior on problems of different sizes.}
\end{figure}

Figure 5 shows a function similar to an increasing exponential.
The proximity of the graph to the origin of coordinates indicates that in problems too
small optimization is hardly effective. This supports the hypothesis
that the results of the heuristic algorithm of insertion is of high quality. In
uncomplicated cases, the initial solution provided by this algorithm hardly
can be improved so that is very close to the optimum.
As the problem grows in size it becomes impossible to obtain 
optimal solutions with heuristic methods using instant execution.
It is therefore effective the use of metaheuristic optimization algorithms
which during determined time track a better solution. How 
we can see in the graph, the simulated annealing algorithm performs an efficient
work on problems of relatively large sizes.

\begin{figure}
\begin{center}
	\includegraphics[height=8cm]{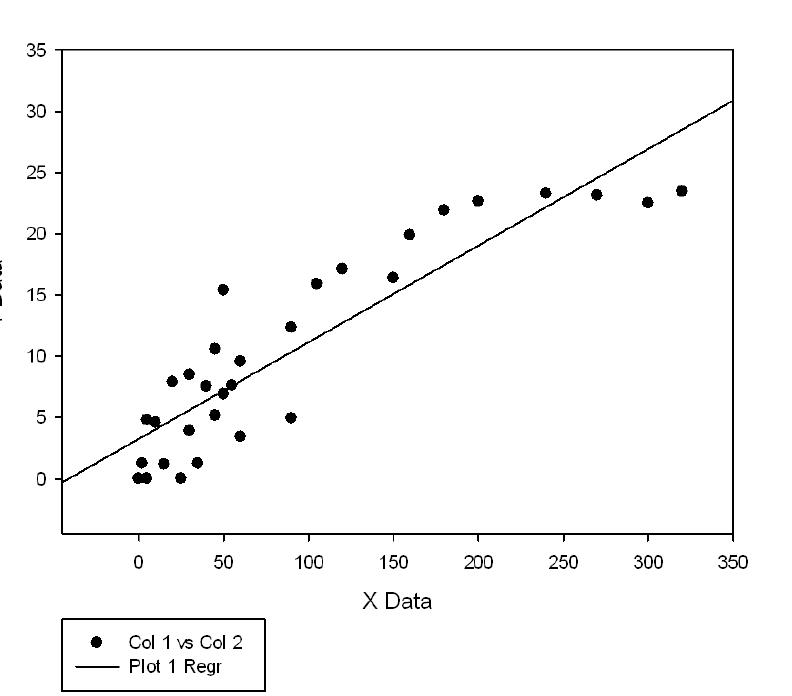}
\end{center}
\caption{Quality results depending on the time needed to generate it.}
\end{figure}

Figure 6 shows in a more specific way the concept of efficiency.
It measures the quality of a result as a function of time necessary to generate it.
It is indicated how it evolves the process of optimization results as
it endows runtime. The vertical axis corresponds to the percentage of
improvement between the initial solution and the optimized solution and the horizontal axis to the
time spent in obtaining the result. In these tests it has been used
always the same problem, of moderate size.

The increasing character of the function tells us how it was assumed that older
runtimes produce higher optimizations. If we deepen the
analysis it highlights the proximity of the graph to the origin of coordinates indicating that
with a very small execution times there are no significant improvements. In those
cases the optimization algorithm has not enough processing for
find a local minimum of the objective function to improve the initial solution.
As the execution time is increased, usually new
local minima that improve the initial solution are emerging and progressively
the identification of best solutions optimize more and more results.

Then the rhythm of
improvement increases and the graph reaches its maximum slope. It is in this position
when the temporary resource is exploited and the appearance of more results
of quality is more frequent.
As the problem is given with large amounts of time, results
improve but at rates increasingly soft. We believe that in these cases, the solution reached is close to the global optimum of the objective function.
The refining process of results can continue but the quality of
the same would not improve ostentatiously.
In view of the obtained results, tests and subsequent
analysis, the initial requirements for that
this application was designed are considered validated and accepted.

\section{Conclusions and final remarks}

We can conclude that our first experimental results show that the
heuristics can be applied effectively with reasonable computational effort. Besides, we
have programmed an interface that allows one to use the algorithms in a friendly way
and that can make possible the actual implementation of our methodology in the cooperative. Furthermore, the comparison between our results and those provided
by the former methodology shows a reduction both in working time and
economic cost. 
Thus, the advantage that can provide the use of this tool is that it provides savings in distribution costs because the distances obtained through an optimization procedure will be lower than those obtained manually without using the program. It also facilitates the task of managing truck routes to the logistics technician who need not spend too much time for it or a complex knowledge of the environment in which she is working. Joint with the distribution problem here studied, there are other aspects associated with possible inefficiencies in the system that might be worth studying (cf. Amiama-Ares et al., 2014 \cite{ami}). One can be the problem of production of feed, which can be another future task to address using tools taken of operations research and computer science.

\section*{Acknowledgements}

This research received financial support from Spanish {\it Ministerio de Econom\'{\i}a y Competitividad} and FEDER through grant MTM2014-53395-C3-2-P that is gratefully acknowledged.

\end{document}